\documentclass{article}

\usepackage{arxiv}

\usepackage[utf8]{inputenc} 
\usepackage[T1]{fontenc}    
\usepackage[]{hyperref}
\usepackage[square, authoryear]{natbib}
\usepackage{url}            
\usepackage{booktabs}       
\usepackage{amsfonts}       
\usepackage{nicefrac}       
\usepackage{microtype}      
\usepackage{amsmath, amsthm, amscd, amssymb, graphicx, color}
\usepackage[labelfont=bf]{caption}
\usepackage{float}
\usepackage{booktabs, tabu, array}
\graphicspath{ {fig/} }

\title{Normality Analysis of Current World Record Computations for Catalan's Constant and Arc Length of a Lemniscate with $a=1$}

\date{August 21, 2019}	

\author{
  Seungmin~Kim\thanks{Corresponding author, ORCID: 0000-0001-8052-723X} \\
  Korea Science Academy of KAIST, 105-47, Baegyanggwanmun-ro, Busanjin-gu, Busan, 47162, Republic of Korea \\
  \texttt{ehf@users.sourceforge.net} \\
}

\begin{document}
\maketitle
\fancypagestyle{empty}{
    \fancyhf{}
    \renewcommand{\headrulewidth}{0pt}
    \cfoot{\thepage}
}
\chead{Kim, S.}

\begin{abstract}
Catalan's constant and the lemniscate constants have been important mathematical constants of interest to the mathematical society, yet various properties are unknown. An important property of significant mathematical constants is whether they are normal numbers. This paper evaluates the normality of decimal and hexadecimal representations of current world record computations of digits for the Catalan's constant (600,000,000,100 decimal digits and 498,289,214,317 hexadecimal digits) and the arc length of a lemniscate with $a=1$ (600,000,000,000 decimal digits and 498,289,214,234 hexadecimal digits). All analyzed frequencies are persistent to the conjecture of Catalan's constant and the arc length of a lemniscate with $a=1$ being a normal number in bases 10 and 16.
\end{abstract}


\textbf{2010 MSC} MSC 11 · MSC 62 · MSC 68

\section{Introduction}
A normal number at a base of $b\in\mathbb Z_{>0}$ is a real number such that each of the b digits is distributed in the equal natural density of \(\frac{1}{b}\). Similarly, for a length $n\in\mathbb Z_{>0}$, all $b^n$ digit combinations in base b have the same natural density of $b^{-n}$, thus being equally likely in occurrence. This implies that no digit or combination of digits occurs more frequently than others. Very few computable mathematical constants have been proved as normal, although the mathematical society widely assumes that computable mathematical constants such as $\pi$, $e$, and $\sqrt{2}$ are normal. 

Catalan's constant and the lemniscate constants have been crucial in the development of mathematics by themselves and their computation methods, and further related discoveries implicate even greater insight. Catalan's constant $G$ is defined as $G=\beta(2)=\sum_{n=0}^{\infty} \frac{(-1)^{n}}{(2n+1)^{2}}$, and whether it is transcendental or even irrational remains an open question in mathematics \citep{nest2016cat}. The arc length of a lemniscate with $a=1$ is defined as $s=\frac{1}{\sqrt{2\pi}}[\Gamma(\frac{1}{4})]^{2}$. The lemniscate constant is defined as $L=\frac{s}{2}=\int_{0}^{1} \frac{dt}{\sqrt{1-t^{4}}}$ and has been proven as a transcendental number \citep{finch2003lemn, todd2000lemn}. The values $L_1=\frac{L}{2}$ and $L_2=\frac{1}{2G}$, where $G$ is Gauss's constant, are sometimes respectively referred to as the first and second lemniscate constants, both values also proven as transcendental \citep{finch2003lemn, todd2000lemn}. 

Catalan's constant, the arc length of a lemniscate with $a=1$, and all lemniscate constants are unknown if it is a normal number. Currently, one of the only plausible methods to approximately verify sequences or irrational numbers are normal is to perform statistical analysis on a large sample of the sequence or number. Because of this, new world records of computed mathematical constants are used for checking statistical consistency for the normality of many important mathematical constants. Trueb \citep{trueb2016digit} has performed an analysis of digit combinations from length one to three on the first $\lfloor{\pi^{e}}\rfloor$ trillion digits of $\pi$ and has verified that the variance of the frequencies overall complies with the expected variance. New world record calculations for the Catalan's constant (600,000,000,100 decimal digits and 498,289,214,317 hexadecimal digits) and the arc length of a lemniscate with $a=1$ (600,000,000,000 decimal digits and 498,289,214,234 hexadecimal digits) as of July 2019 have recently been calculated and verified using the \textit{y-cruncher} program \citep{yeecruncher, kim2019cat, kim2019lemn}. This study attempts to expand the insight of these two important mathematical constants on the conjecture of the Catalan's constant and the arc length of a lemniscate with $a=1$ by statistically analyzing the calculated world record digit computation results \citep{kim2019cat, kim2019lemn}. 

\section{Materials and Methods}
Analysis of digits was done using a custom-coded Python 3 script, speed optimized using the Portable PyPy3.6 v7.1.1 JIT compiler \citep{pypy36port} on CentOS 7.4 using an Intel Xeon (Skylake Purley) CPU. The Python script was coded so that the integrity of arbitrary length digit combinations is ensured if single digit counts are correct, and the single-digit counts were initially cross-checked for verification using the Digit Viewer application of \textit{y-cruncher} \citep{yeecruncher}. Statistical analysis and visualization were done using a modified version of the method and code used in Trueb \citep{trueb2016digit} using the CERN ROOT Toolkit \citep{root1997} Release 5.34/38, compiled using GCC version 7.4.0 on Ubuntu 18.04 LTS using a KVM virtual machine with an Intel Xeon (Haswell) CPU. The source code used for data analysis was based on the ANSI C99 standard of the C language. 

Digit occurrences from length one to three were counted starting from the decimal point. The decimal and hexadecimal expressions of the Catalan's constant and the arc length of a lemniscate with $a=1$ were respectively used for the analysis of digit occurrences. The computation results for Catalan's constant had 600,000,000,100 decimal digits and 498,289,214,317 hexadecimal digits for each base expression after the decimal point, and the computation results for the arc length of a lemniscate with $a=1$ had 600,000,000,000 decimal digits and 498,289,214,234 hexadecimal digits for each base expression after the decimal point. The expected variance of the frequencies and the expected error of the variance has been calculated by assuming the frequencies followed a binomial distribution around the limiting frequency of $b^{-k}$ \citep{trueb2016digit}. Accessory regions were plotted with the area between the two vertical lines closest to the center of the figure representing the region within one standard deviation and the remaining area between two remaining vertical lines representing the region between one and two standard deviations \citep{trueb2016digit}. 

\section{Results}
Figures 1 to 6 depict the frequency distributions of all sequences from a length of one to three for the decimal and hexadecimal representations of Catalan's constant in the form of a histogram. Table 1 lists the predicted and actual variances for the frequency distributions of Catalan's constant. 

Figures 7 to 12 depict the frequency distributions of all sequences from a length of one to three for the decimal and hexadecimal representations of the arc length of a lemniscate with $a=1$ in the form of a histogram. Table 2 lists the predicted and actual variances for the frequency distributions of the arc length of a lemniscate with $a=1$. 

\begin{figure}[H]
    \centering
    \begin{minipage}[b]{0.45\textwidth}
        \includegraphics[width=\textwidth]{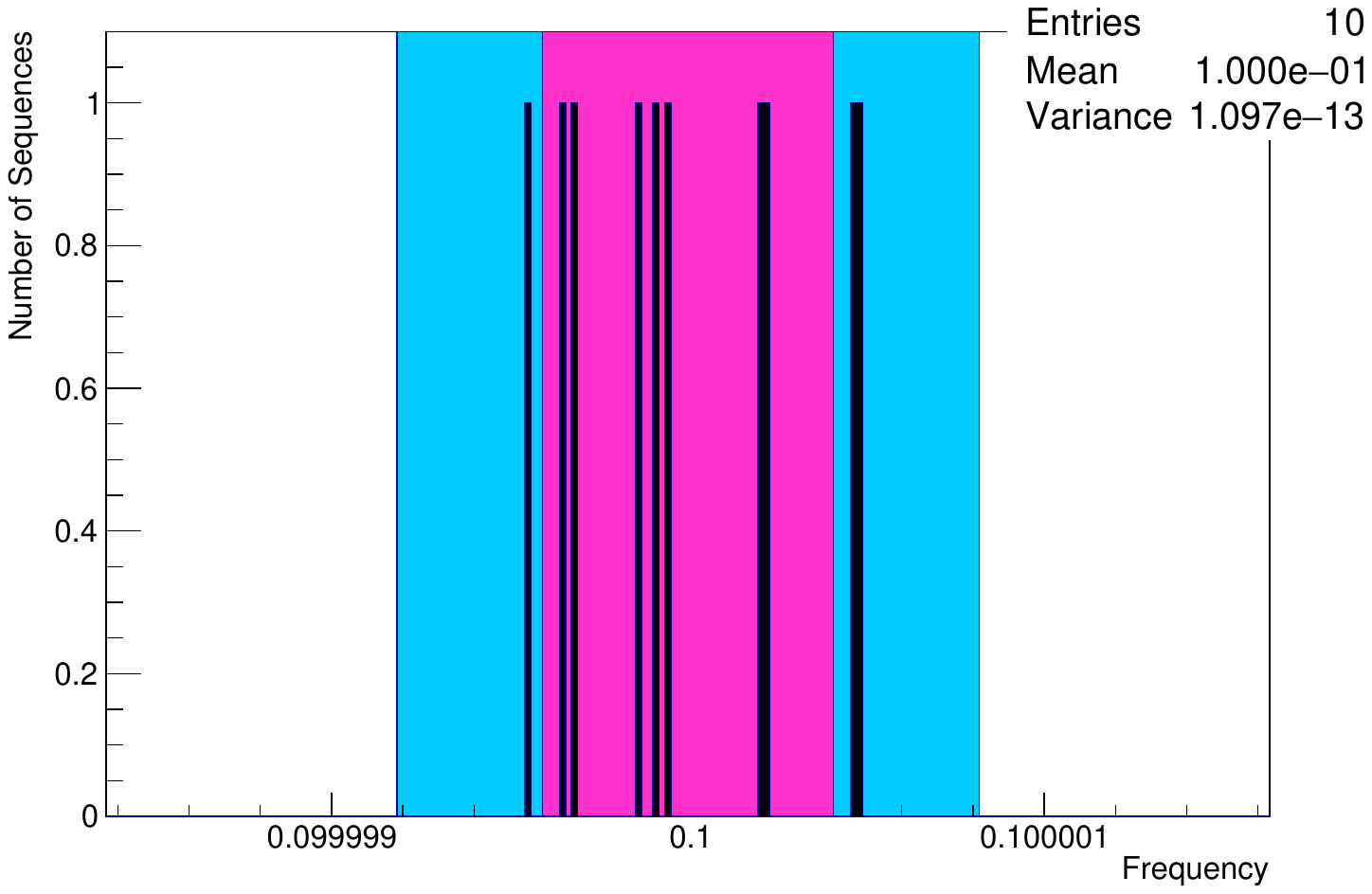}
        \caption{Frequencies of all combinations of length 1 (digits 0--9) in the decimal representation of the Catalan's constant.}
        \label{fig:1}
    \end{minipage}\hfill
    \begin{minipage}[b]{0.45\textwidth}
        \includegraphics[width=\textwidth]{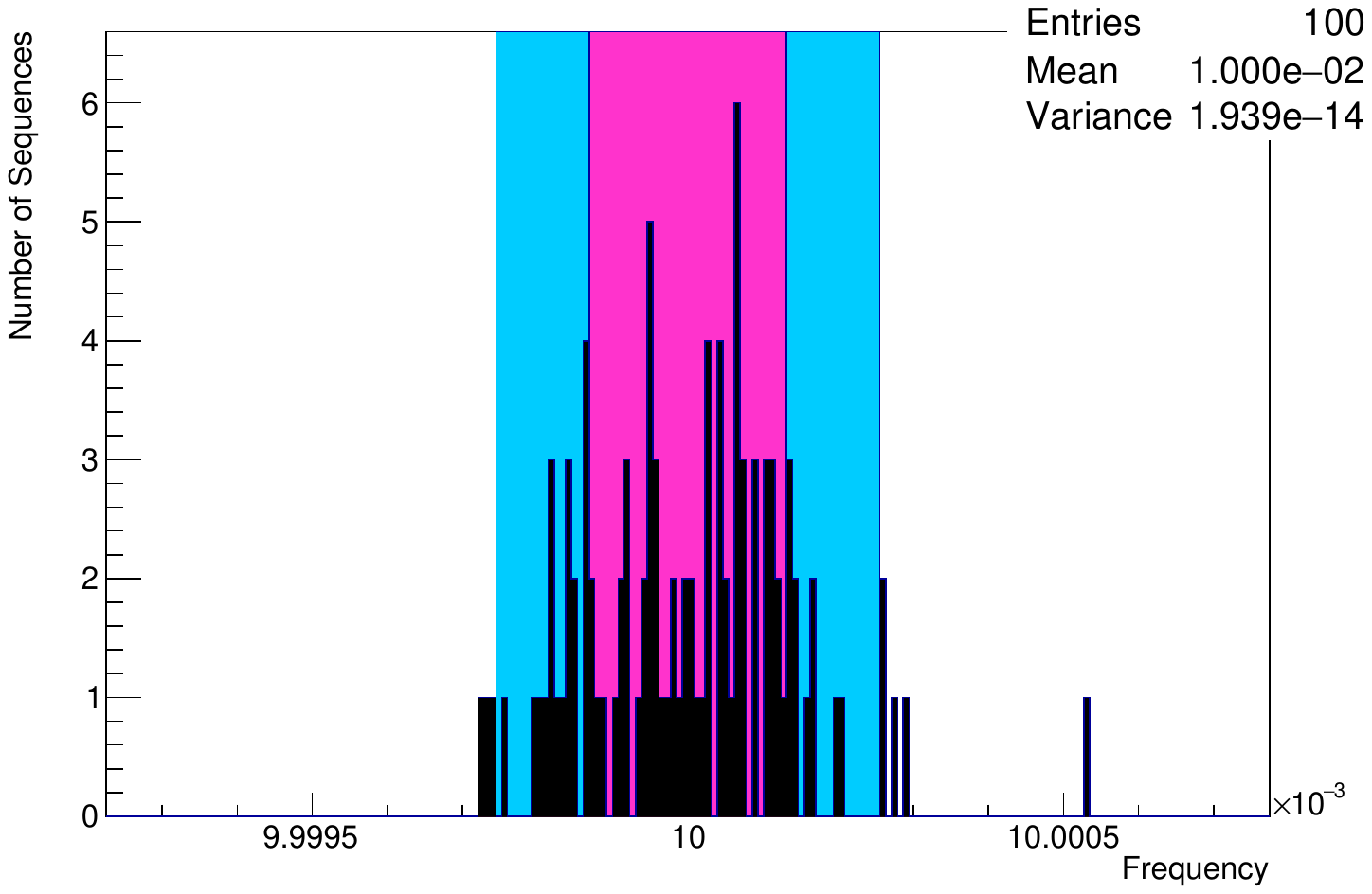}
        \caption{Frequencies of all combinations of length 2 (00--99) in the decimal representation of the Catalan's constant.}
        \label{fig:2}
    \end{minipage}
\end{figure}

\begin{figure}[H]
    \centering
    \begin{minipage}[b]{0.45\textwidth}
        \includegraphics[width=\textwidth]{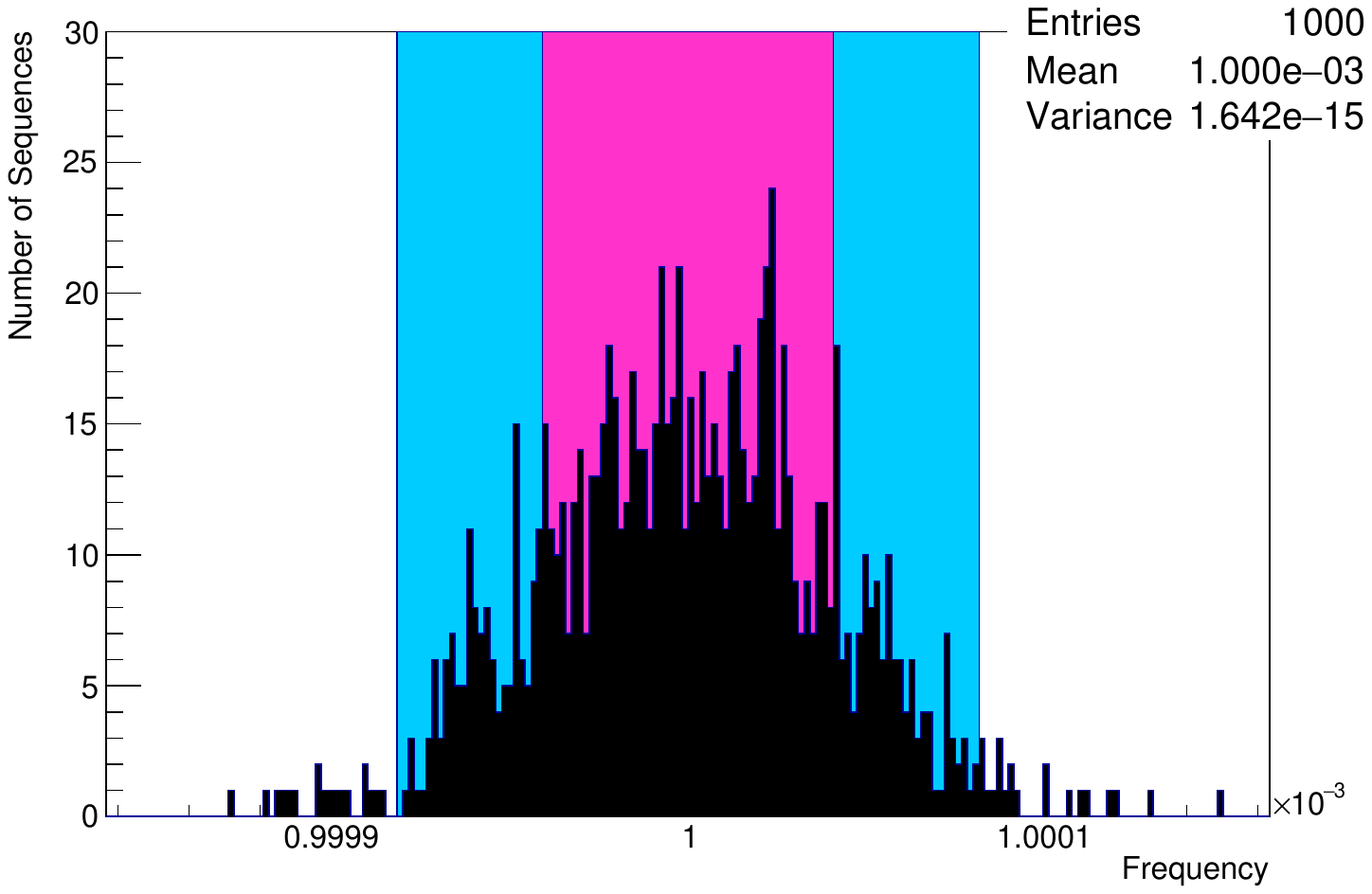}
        \caption{Frequencies of all combinations of length 3 (000--999) in the decimal representation of the Catalan's constant.}
        \label{fig:3}
    \end{minipage}\hfill
    \begin{minipage}[b]{0.45\textwidth}
        \includegraphics[width=\textwidth]{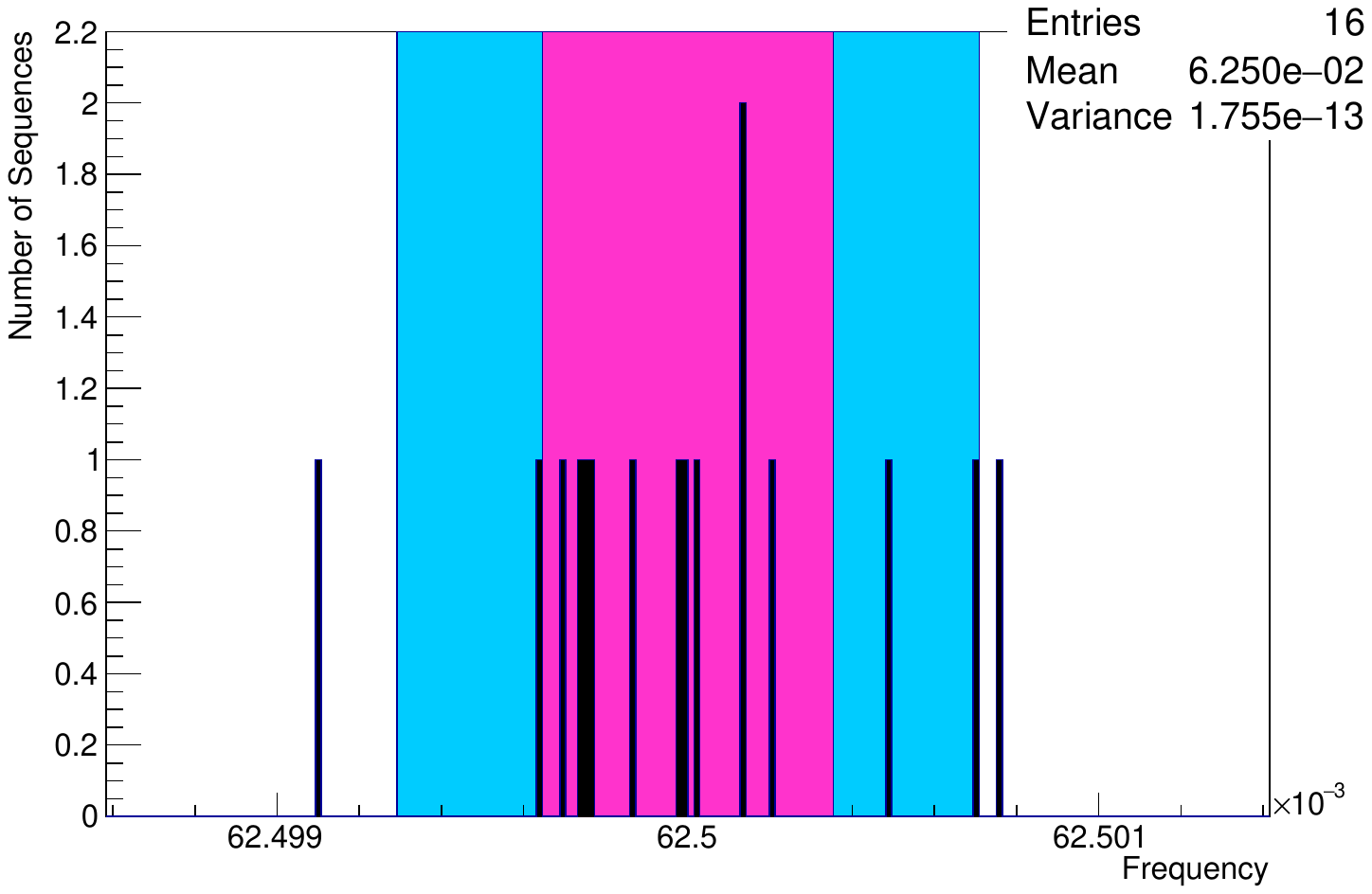}
        \caption{Frequencies of all combinations of length 1 (digits 0--F) in the hexadecimal representation of the Catalan's constant.}
        \label{fig:4}
    \end{minipage}
\end{figure}

\begin{figure}[H]
    \centering
    \begin{minipage}[b]{0.45\textwidth}
        \includegraphics[width=\textwidth]{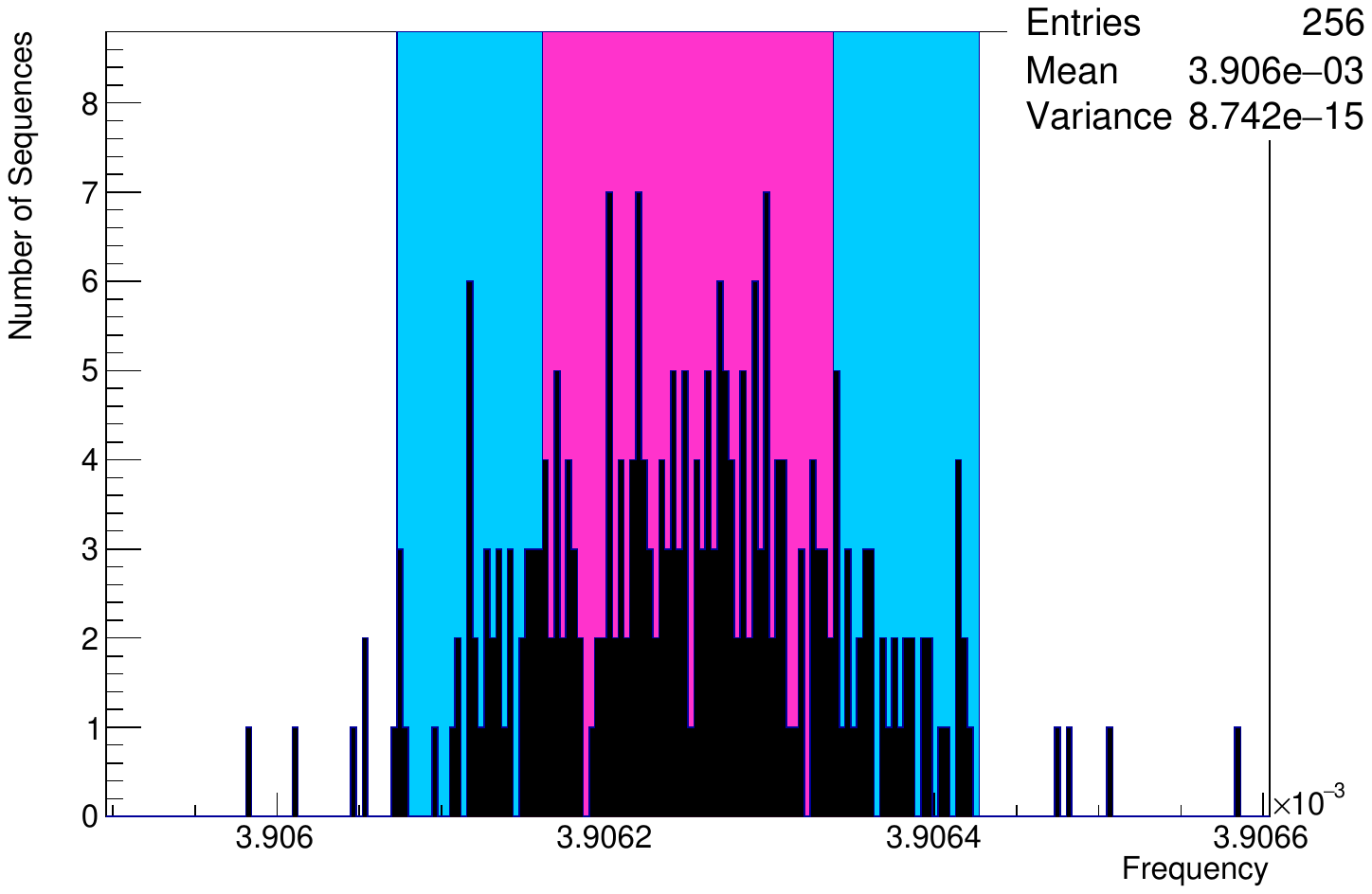}
        \caption{Frequencies of all combinations of length 2 (00--FF) in the hexadecimal representation of the Catalan's constant.}
        \label{fig:5}
    \end{minipage}\hfill
    \begin{minipage}[b]{0.45\textwidth}
        \includegraphics[width=\textwidth]{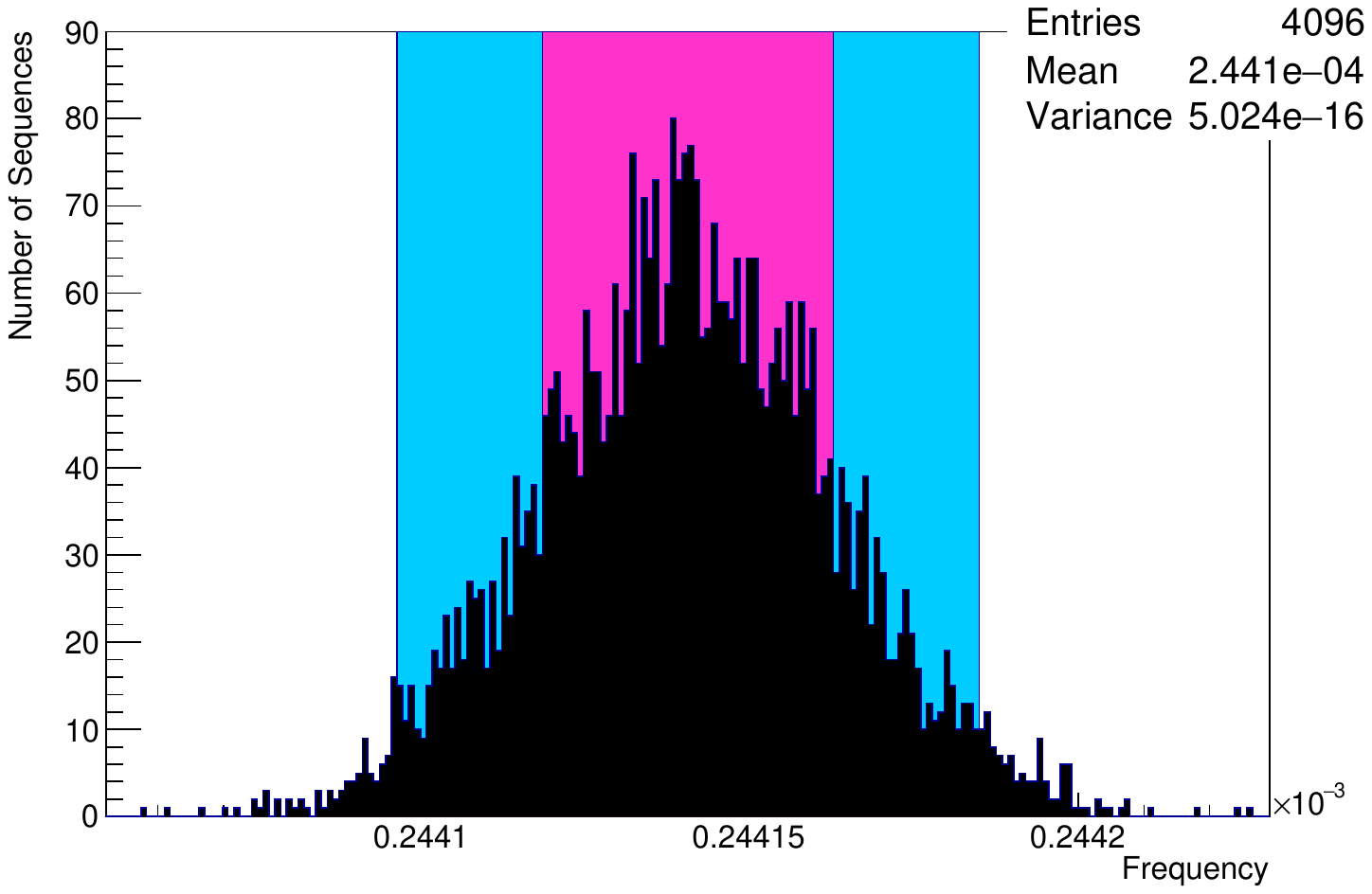}
        \caption{Frequencies of all combinations of length 3 (000--FFF) in the hexadecimal representation of the Catalan's constant.}
        \label{fig:6}
    \end{minipage}
\end{figure}

\begin{table}[H]
\centering
\caption{Predicted and actual variances of frequencies of all sequences of length 1--3 in the decimal and hexadecimal representations of Catalan's constant.}
\begin{tabu} to \textwidth{@{}|c|c|X[c]|X[c]|c|@{}}
\toprule
    Base & Length of Sequence & Predicted Variance and Error of Frequencies & Actual Variance of Frequencies & Deviation [$\sigma$] \\ \midrule
    10 & 1 & $(1.500\pm0.707)\times10^{-13}$ & $1.097\times10^{-13}$ & $0.570$ \\ \midrule
    10 & 2 & $(1.650\pm0.235)\times10^{-14}$ & $1.939\times10^{-14}$ & $-1.232$ \\ \midrule
    10 & 3 & $(1.665\pm0.074)\times10^{-15}$ & $1.642\times10^{-15}$ & $0.306$ \\ \midrule
    16 & 1 & $(1.176\pm0.429)\times10^{-13}$ & $1.755\times10^{-13}$ & $-1.349$ \\ \midrule
    16 & 2 & $(7.809\pm0.692)\times10^{-15}$ & $8.742\times10^{-15}$ & $-1.349$ \\ \midrule
    16 & 3 & $(4.898\pm0.108)\times10^{-16}$ & $5.024\times10^{-16}$ & $-1.160$ \\ \midrule
    \bottomrule
\end{tabu}
\label{tab:1}
\end{table}

\begin{figure}[H]
    \centering
    \begin{minipage}[b]{0.45\textwidth}
        \includegraphics[width=\textwidth]{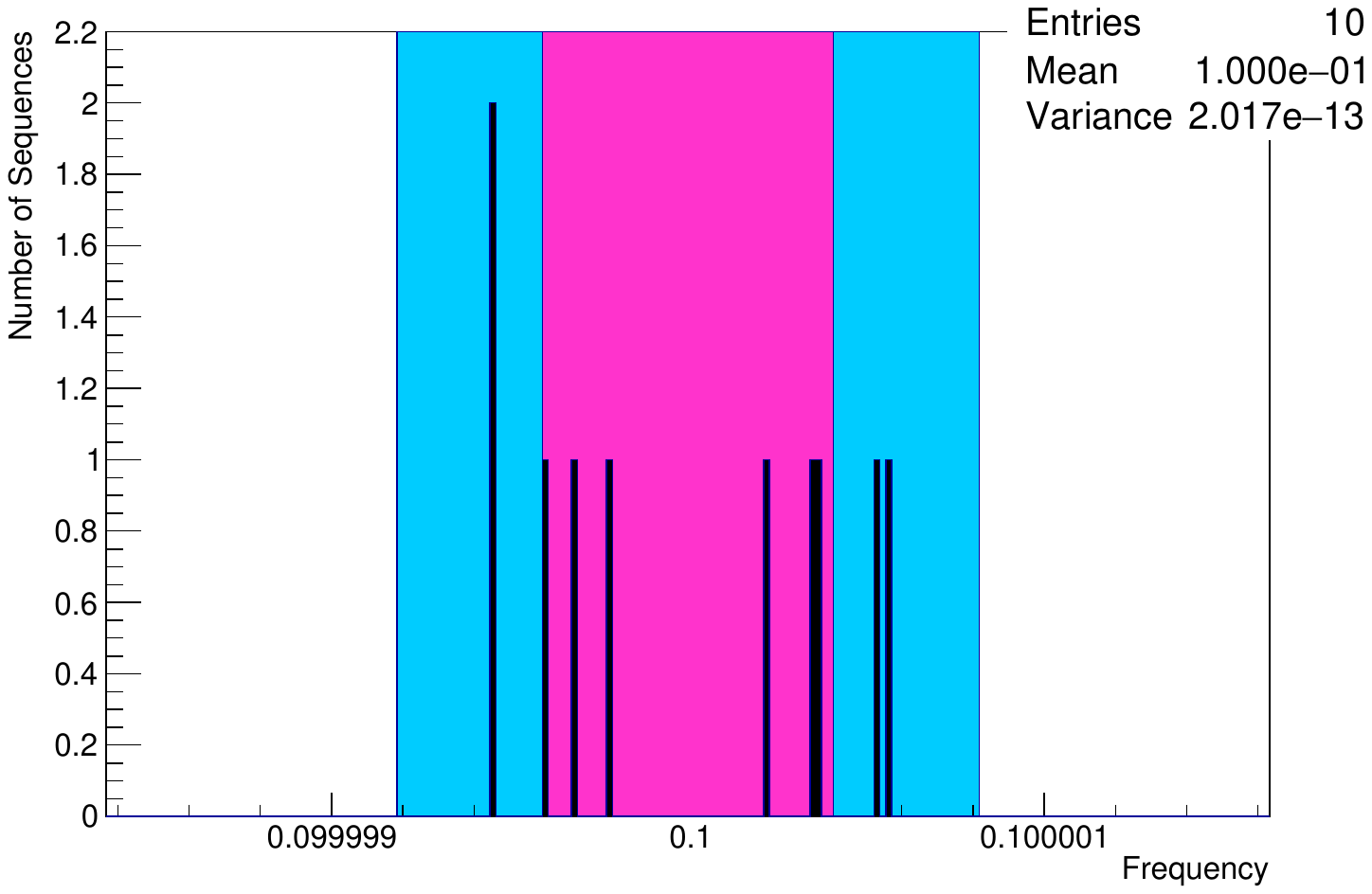}
        \caption{Frequencies of all combinations of length 1 (digits 0--9) in the decimal representation of the arc length of a lemniscate with $a=1$.}
        \label{fig:7}
    \end{minipage}\hfill
    \begin{minipage}[b]{0.45\textwidth}
        \includegraphics[width=\textwidth]{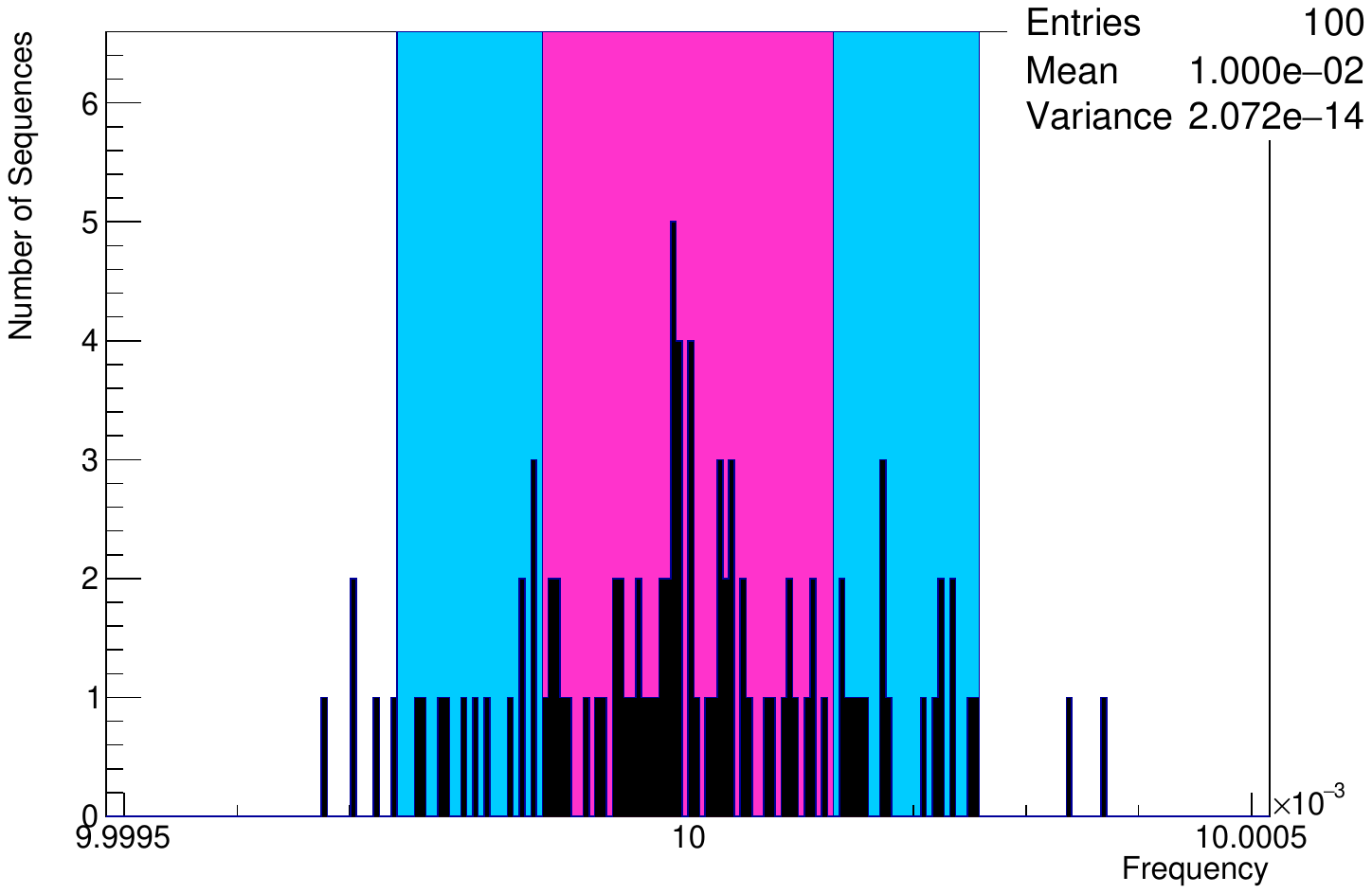}
        \caption{Frequencies of all combinations of length 2 (00--99) in the decimal representation of the arc length of a lemniscate with $a=1$.}
        \label{fig:8}
    \end{minipage}
    \begin{minipage}[b]{0.45\textwidth}
        \includegraphics[width=\textwidth]{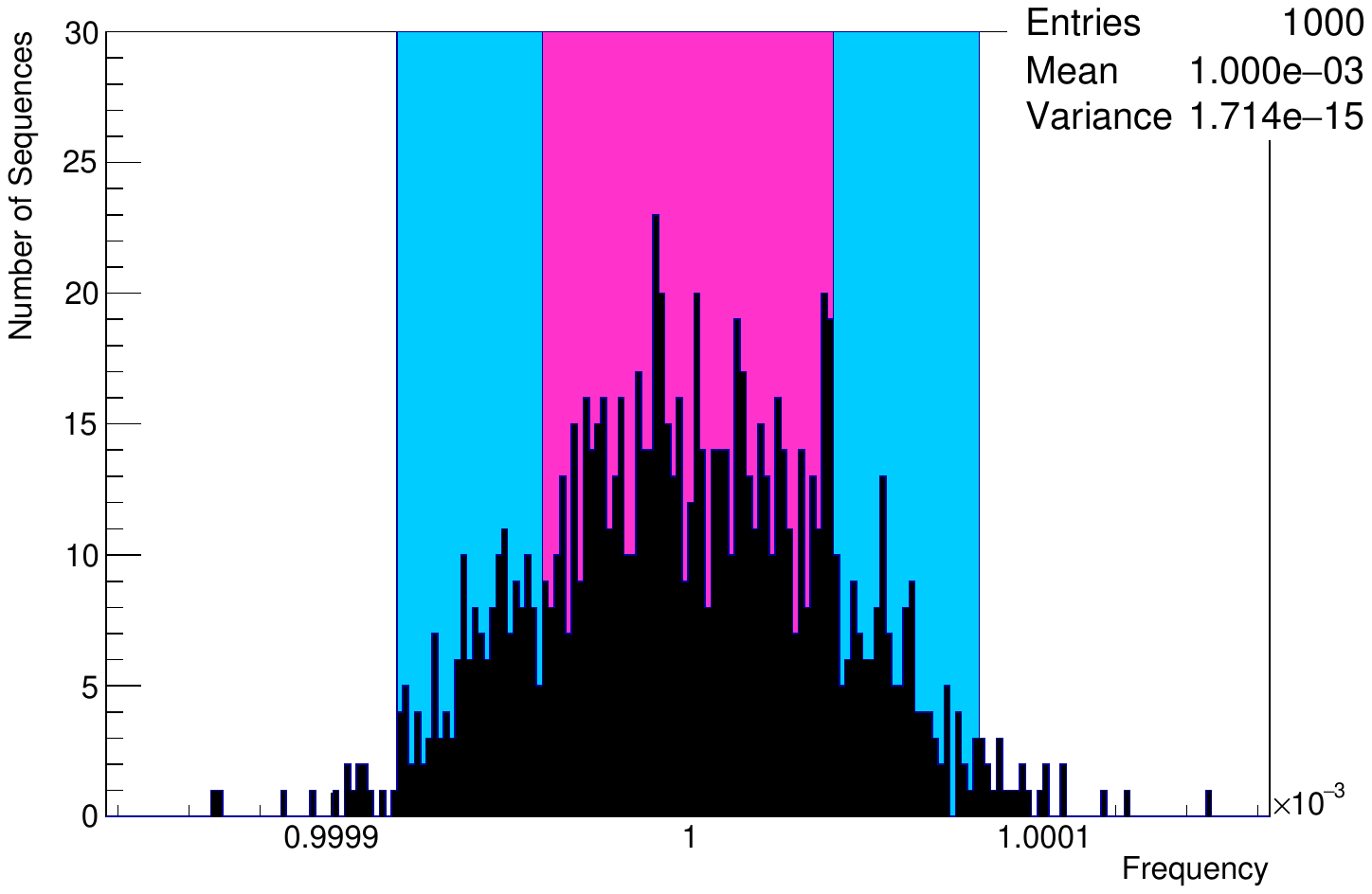}
        \caption{Frequencies of all combinations of length 3 (000--999) in the decimal representation of the arc length of a lemniscate with $a=1$.}
        \label{fig:9}
    \end{minipage}\hfill
    \begin{minipage}[b]{0.45\textwidth}
        \includegraphics[width=\textwidth]{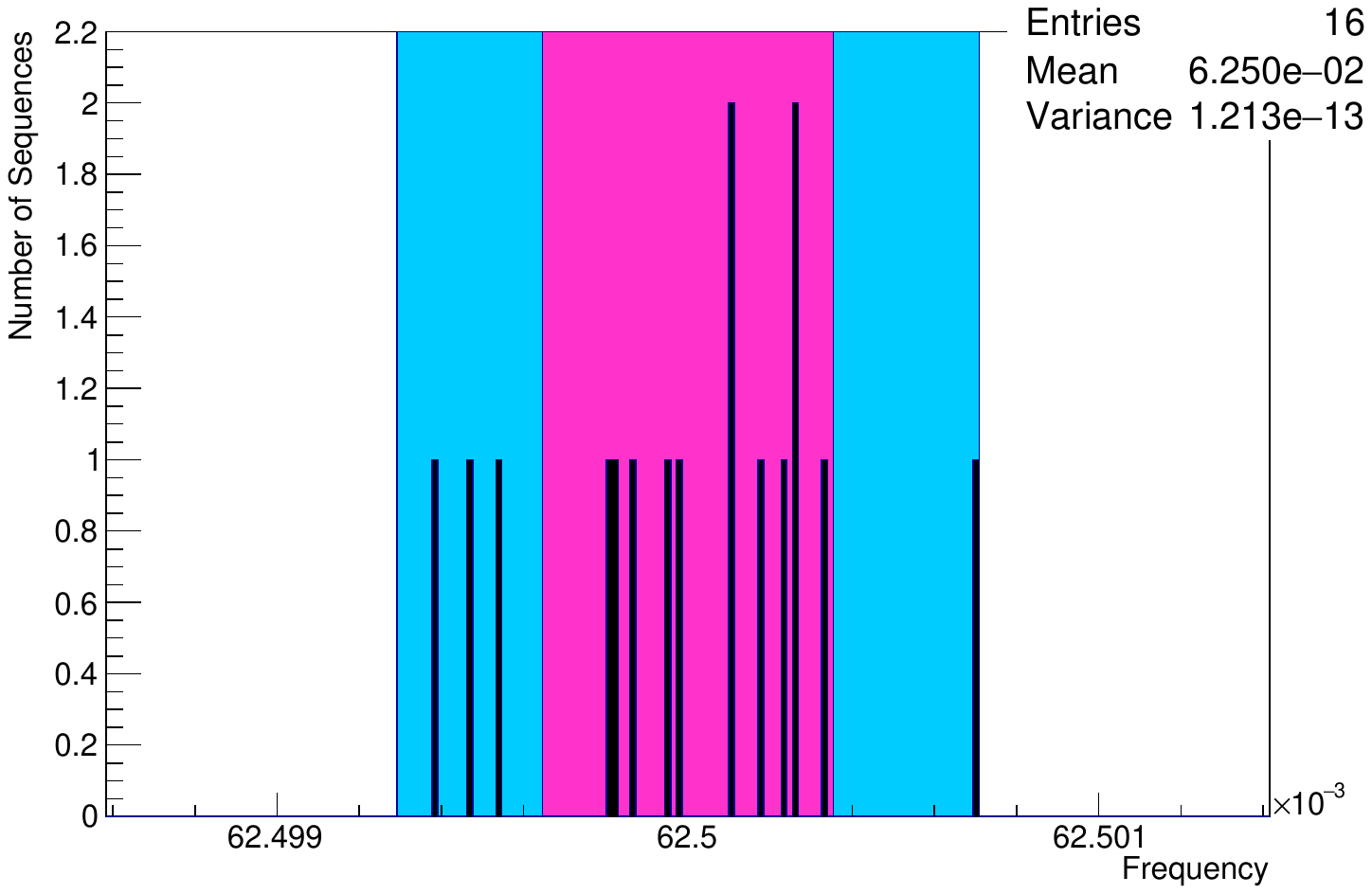}
        \caption{Frequencies of all combinations of length 1 (digits 0--9) in the hexadecimal representation of the arc length of a lemniscate with $a=1$.}
        \label{fig:10}
    \end{minipage}
    \begin{minipage}[b]{0.45\textwidth}
        \includegraphics[width=\textwidth]{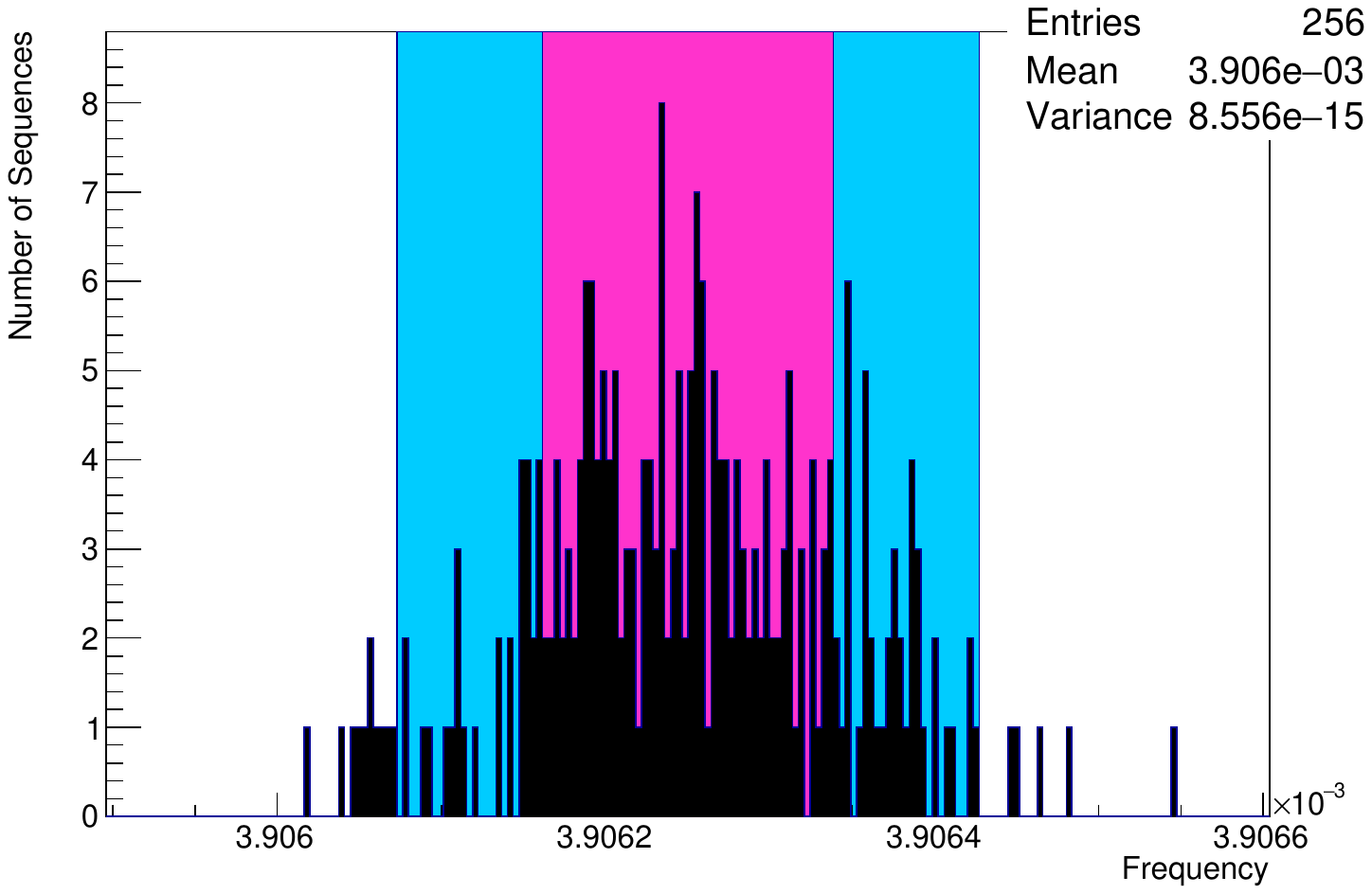}
        \caption{Frequencies of all combinations of length 2 (00--99) in the hexadecimal representation of the arc length of a lemniscate with $a=1$.}
        \label{fig:11}
    \end{minipage}\hfill
    \begin{minipage}[b]{0.45\textwidth}
        \includegraphics[width=\textwidth]{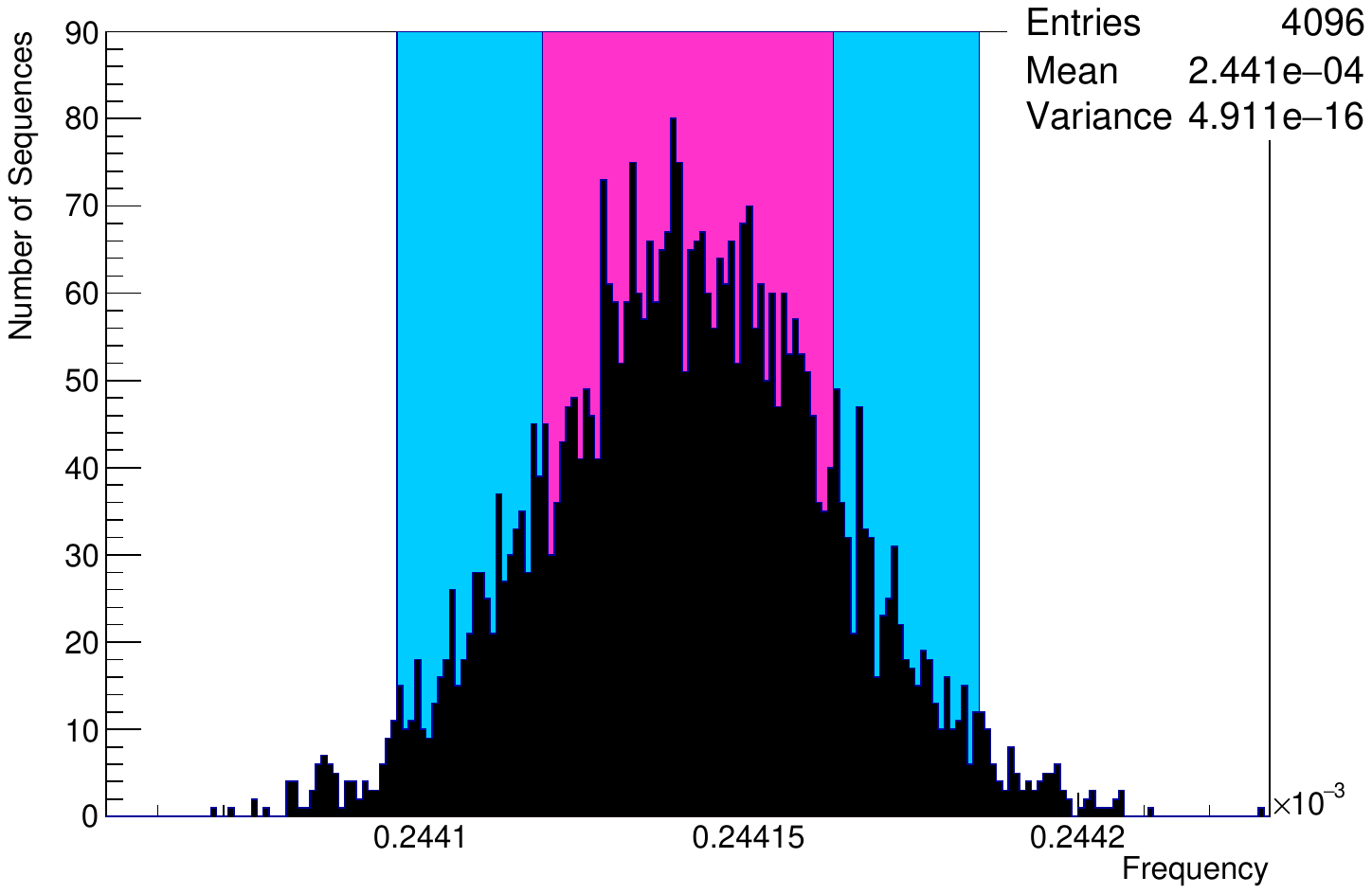}
        \caption{Frequencies of all combinations of length 3 (000--999) in the hexadecimal representation of the arc length of a lemniscate with $a=1$.}
        \label{fig:12}
    \end{minipage}
\end{figure}

\begin{table}[H]
\centering
\caption{Predicted and actual variances of frequencies of all sequences of length 1--3 in the decimal and hexadecimal representations of the arc length of a lemniscate with $a=1$.}
\begin{tabu} to \textwidth{@{}|c|c|X[c]|X[c]|c|@{}}
\toprule
    Base & Length of Sequence & Predicted Variance and Error of Frequencies & Actual Variance of Frequencies & Deviation [$\sigma$] \\ \midrule
    10 & 1 & $(1.500\pm0.707)\times10^{-13}$ & $2.017\times10^{-13}$ & $-0.731$ \\ \midrule
    10 & 2 & $(1.650\pm0.235)\times10^{-14}$ & $2.072\times10^{-14}$ & $-1.798$ \\ \midrule
    10 & 3 & $(1.665\pm0.074)\times10^{-15}$ & $1.714\times10^{-15}$ & $-0.660$ \\ \midrule
    16 & 1 & $(1.176\pm0.429)\times10^{-13}$ & $1.213\times10^{-13}$ & $-0.086$ \\ \midrule
    16 & 2 & $(7.809\pm0.692)\times10^{-15}$ & $8.556\times10^{-15}$ & $-1.080$ \\ \midrule
    16 & 3 & $(4.898\pm0.108)\times10^{-16}$ & $4.911\times10^{-16}$ & $-0.113$ \\ \midrule
    \bottomrule
\end{tabu}
\label{tab:2}
\end{table}

\section{Discussion}
The frequencies of sequences from length 1 to 3 for both decimal and hexadecimal representations in both Catalan's constant and the arc length of a lemniscate with $a=1$ overall coincide with the expected distribution. The distribution of frequencies for Catalan's constant has a maximum deviation of $1.232$ standard deviations (absolute value) for the decimal representation, and $1.349$ standard deviations (absolute value) for the hexadecimal representation. The distribution of frequencies for the arc length of a lemniscate with $a=1$ has a maximum deviation of $1.798$ standard deviations (absolute value) for the decimal representation, and $1.080$ standard deviations (absolute value) for the hexadecimal representation. 

These results for a very large sample data of around 600 billion digits of both Catalan's constant and the arc length of a lemniscate with $a=1$ are statistically persistent to the conjecture of Catalan's constant and the arc length of a lemniscate with $a=1$ being a normal number in bases 10 and 16. 

In this study, it has been possible to make a statistical assumption of whether Catalan's constant and the arc length of a lemniscate with $a=1$ is a normal number. Further approaches based on the methods used in this study can enable the mathematical society to gain evidence on the normality of a large range of computable mathematical constants of interest. 

\section*{Acknowledgment}
The author greatly thanks Dr. Peter Trueb for providing code and insights to formulating the methods for data analysis in this study, Mr. Alexander J. Yee for providing insights related to his program \textit{y-cruncher}, and Dr. Ian Cutress for verifying the world record computation for the arc length of a lemniscate with $a=1$. 

\section*{Data Availability Statement}
The data that support the findings of this study are openly available in the Internet Archive at \url{https://archive.org/details/catalan_190618}, reference number \texttt{ark:/13960/t9f55d078}, and at \url{https://archive.org/details/lemworldrec_190512}, reference number \texttt{ark:/13960/t56f3gj32}.

\bibliographystyle{unsrtnat}  
\bibliography{references}  



\end{document}